# Randomness of D Sequences via Diehard Testing

James Bellamy

**Abstract.** This paper presents a comparison of the quality of randomness of D sequences based on diehard tests. Since D sequences can model any random sequence, this comparison is of value beyond this specific class.

**Introduction**
Randomness of binary sequences has been studied for a long time [1] and it continues to interest researchers because of the need to distinguish (in applicable situations) between randomness that has a classical basis from randomness that emerges out of quantum behavior. Random binary sequences are important not only in classical cryptography but also in quantum cryptography since the choice of random bases in BB84 [2] or polarization rotation in the three-stage protocol [3] must be done classically. One class of random sequences that can model a wide variety of random sequences is that of D sequences [4]-[6]. Here, we examine the randomness of a certain class of D sequences.

**D Sequence Algorithm**
D sequences are capable of generating random numbers, as outlined by [7]. The class of D sequences that will be tested in this paper is described as follows. First, a divisor $m$ will be chosen. Next, choose a value $x$ which will be the multiplier and base value. The algorithm will iterate along computing $a(i)$ as such:
$$a(i) = (l^i \bmod m)$$
This sequence of values can then be converted to a binary stream $b(i)$ by:
$$b(i) = (a(i) \bmod m)$$
This is a general description of a D Sequence. The first step I completed was implementing the modified D Sequence Algorithm in MATLAB. This was a fairly straightforward procedure. I ran into one major issue with this implementation but this also gave me the idea to make an adjustment to the algorithm that actually led to a faster implementation.

The modified algorithm first takes a maximum range which can be inputted by the user and finds the two highest prime numbers that fall into that range. It then multiplies these two values in order to get the value for $m$. Next, it chooses a value that is approximately 1/4 of the modulus that will be used. After it is completely set up, it will increment through all values which are in this range and compute all $a(i)$ of the unique values until there is a repeat. Once it has found all of the unique values it will drop out of the loop and move into the next step, which simplifies all of the values down to their LSB and add these to the next element.

**Diehard Tests**
Diehard tests consist of a battery of 15 different individual tests, which are used to evaluate how well a pseudo-random number generator produces values. They were originally developed by George Marsaglia and published in 1995. Here I outline a modified implementation of the D

Sequence Algorithm for pseudo-random number generation and the analysis of the quality of the numbers generated using the diehard tests.  This will then be compared to two different pseudo-random number generators and a discussion of the results will be conducted.

In order to run the diehard tests we need to have a file containing 11-12 megabytes worth of pseudo-random number values.  Therefore, the number of elements that has to be generated in the first loop is actually 8 times that value, based on the fact that the algorithm currently only takes one bit from each value generated.  Using MATLAB's built in speed test; this takes upward of 15 minutes to create the entire file, as shown in Fig. 1.

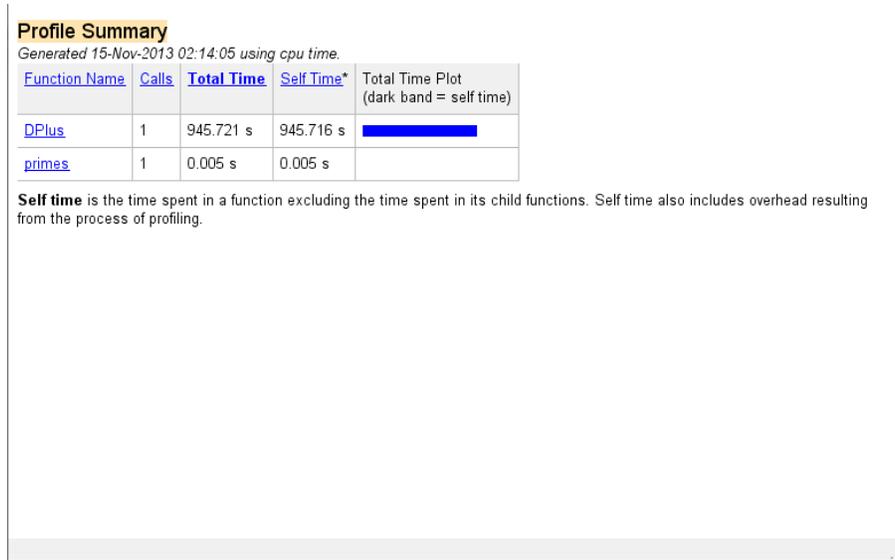

Figure 1

Now that the input file has been created, we are able to move onto actually performing the diehard tests on the random data set that we have just created.  The Diehard test battery consists of 15 different tests, each of which tests the data set in a different way.  These tests are as follows, and are described as such by [9] and [10]:

**Birthday Spacings Test:**   Chooses *m* birthdays in a year of *n* days, listing the spacings between the birthdays.  It then determines a best-fit line and determines the deviation from the expected distribution function.  This deviation determines the quality of this test.

**Overlapping 5-Permutation Test:** Looks at a sequence of one million random integers.  From this, it determines for each set of 5 consecutive integers what the ordering of the five is, (in ascending or descending order.)  It then counts the number of each of the 120 (5!) different configurations and finds the variance from the expected plot of these values.

**Binary Rank Test**: 40,000 31x31 matrices are made from 31 different integers and then the rank of these matrices is found.  The count of matrices at rank 28, 29, 30, and 31 is then found and compared to expected values.  This is then performed 100,000 times for 6x8 matrices through 25 iterations, each returning a p-value.

**Bitstream Test:** Counts the number of missing 20-bit words in a 2^21 bit bitstream. This should be normally distributed.

**OPSO:** Similar to bitstream test, but uses 10-bit letters and 2-letter words and takes different bits from a word each time.

**OQSO (Overlapping-Quadruples-Sparse-Occupancy):** Again similar, but uses 5-bit letters and 4-letter words. The p-value is calculated based on variance from a mean.

**DNA Test:** Again this is similar to the last three, but it uses 5-bit letters and 4-letter words.

**Count the 1's Test:** Uses a sample of 256,000, and determines how many times 8-bit letter and 5-letter words appear and calculates variance. There are two of these tests: one on a stream of bytes and the other on specific bytes.

**Parking Lot Test:** Randomly park cars in a square of side 100, and if a 'crash' occurs, do not park the car. Continue this for $n$=12,000 cars and compare how many successfully parked to the expected result.

**Minimum Distance Test:** Chooses $n$=8,000 random points in a square and then finds $d$, which is the distance between all pairs of points. These distances should be exponentially distributed. P-value is given by the variance from the expected value.

**3D Spheres Test:** Choose 4000 random points placed in a 1000x1000 cube. At each point center a sphere large enough to just touch the nearest cube. The minimum radius is found and this test is repeated 20 times. These radii should be exponentially distributed and again the p-value will be determined based on how far off the results are from the expected.

**Squeeze Test:** Starting with $k$ = 2^31, multiple by floating point values on [0, 1) until $k$ = 1. This is repeated 100,000 times and the number of times it is less than 7 and greater than 47 is counted and compared to an expected exponential to get the p-value.

**Overlapping Sums Test**: Add 100 consecutive integers, incrementing by one element at a time and storing each of the sums, which should have a normal distribution. Compare the results to a normal distribution to get p-value.

**Runs Test:** Iterates through 10,000 values and determines consecutive 'runs' where elements are increasing in value or decreasing in value consecutively. This is then compared to a well-known covariance matrix.

**Craps Test:** Plays 200,000 games of craps, finding the number of throws necessary to end each game. Number of wins should be very close to a normal distribution.

Essentially, the desire is for the p-values returned by each of the tests to have a normal distribution on (0, 1). Therefore, it is highly unlikely for a set of random numbers to generate a large percentage of values which are extremely close to either 0 or 1. However, there is a disclaimer at the top of the documentation that points out that values for which $p < 0.025$ or $p > 0.975$ do not necessarily indicate failure when they are isolated to a single case. It is only when there are a higher number of extreme *p*-values which are within 5-6 decimal places of the edges that you should consider the test a failure. [8], [11]

In addition to testing the modified D Sequence Algorithm which I had written myself, I also tested a couple of other algorithms which were already implemented and made claims to pass the diehard tests. These two functions are the Mersenne Twister and the KISS algorithms.

The first additional PRNG that I decided to use was the Mersenne Twister algorithm. This is the standard pseudo-random number generator used by many different programming languages, including MATLAB, PHP, and Python. It is named after the fact that the period is equal to a Mersenne Prime. I decided to use MATLAB's implementation of this function, and in their implementation the Mersenne Prime is is $2^{19937} - 1$. [12], [13], [14]

The second PRNG I used was the KISS (Keep It Simple, Stupid) Algorithm, which is implemented by George Marsaglia. As the name implies, this algorithm is very simple. It combines two 16-bit Multiply With Carry (MWC) generators with a 3-Shift Register generator (3SRG) and a Congruential generator. The MWCs are combined and then XORed with the Congruential Generator which is then added to the 3SRG. Basically, this generator combines three extremely basic PRNG algorithms in order to create an algorithm that overcomes all of the deficiencies of the other three: It increases the period of the MWC from $2^{60}$ to $2^{123}$, It has much better variability among the lower 16 bits than the Congruential, and it avoids the linear dependence between bits from the 3SRG. [15]

In the following I will discuss all of the tests which failed in some capacity and attempt to give some thoughts. I felt it unnecessary to discuss on a test by test basis the fact that all three PRNGs passed a given test, as that was the case on the majority of the occasions. For some tests, a single p-value is returned, and whether the test was passed is fairly cut and dry. Others return dozens of p-values, so we need to look at the values as a whole and determine if the distribution is as expected.

*Overlapping 5-Permuation Test:*

| KISS | D Sequence | Twister |
|---|---|---|
| .079254 | .248004 | .510352 |
| .898063 | .004112 | .475801 |

Figure 2

If I were to set the threshold at .005<p<.095, then as you can see from the highlighted value in Fig. 2, the D Sequence will fail this test. However, based on the description given by Marsaglia himself, this alone is not enough to indicate a failure of the test. It is outside of the typical bounds though and should be noted in the event that more tests fail. In the event that more than one test fails for the same PRNG and only that PRNG, I prefer to consider the frequency of failure of a given test and increase the weight given to what otherwise may have been considered a minor failure.

*Binary Rank Test*

First, I plotted a histogram and empirical cumulative distribution function of the 25 p-values returned by the 6x8 binary rank test for each of the PRNGs, and the results are shown in Fig. 3:

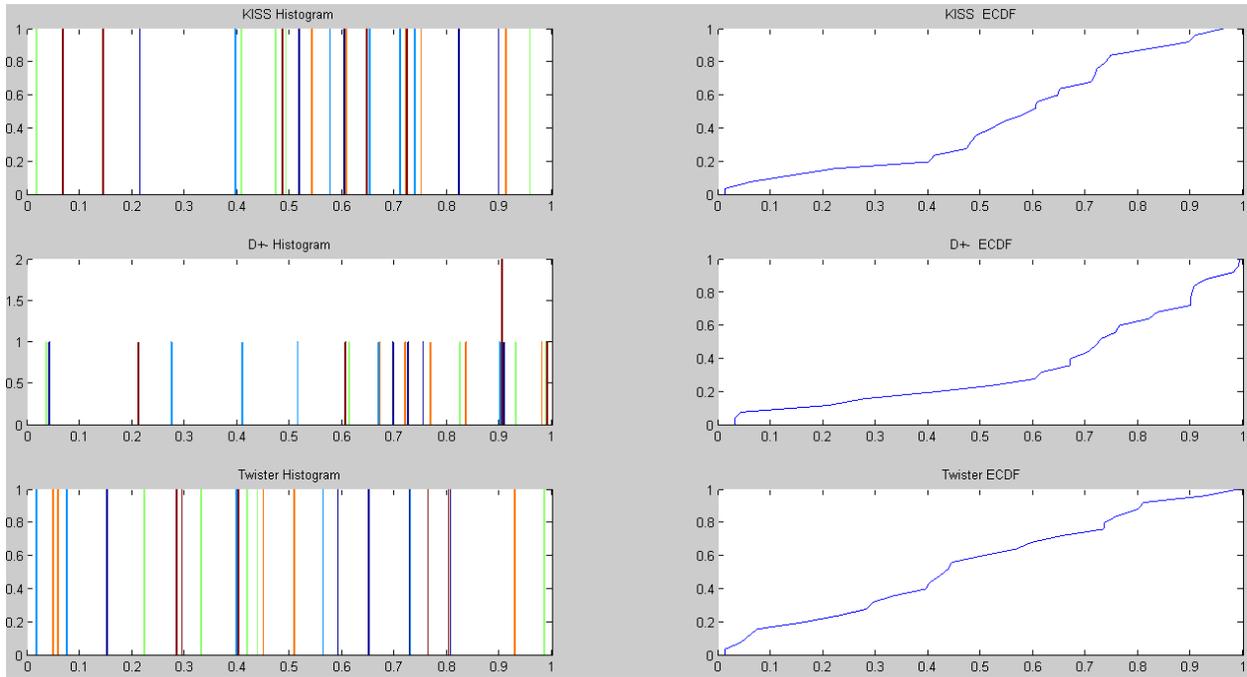

The main difference I notice is that the p-values for the D Sequence appear to be concentrated slightly higher than would be expected with a uniform distribution. The p-values that are calculated from these distributions are shown in the overall 6x8 row in Fig. 4 on the following page:

|  | KISS | D Sequence | Twister |
|---|---|---|---|
| 31x31 | .671639 | .396649 | .997521 |
| 32x32 | .888721 | .322291 | .366271 |
| 6x8 overall | .757647 | .997771 | .141930 |

Figure 4

This time, for two different PRNGs one of the three p-values was outside of the desirable range. There are two directions that this leads my thinking on these results. First of all, the fact that the Mersenne Twister, which as mentioned earlier is becoming more and more dominant in programming languages as the default generator is of note. To me the fact that this PRNG that is held in such high esteem fell well outside of the desired range strengthens my belief that a rare failure of a test should be treated as a false positive until a larger number of failures proves otherwise. However, the fact that the D Sequence algorithm has failed "twice" before the KISS algorithm has "failed" once also leads me to further question the viability of the D Sequence algorithm.

*Bitstream Test*

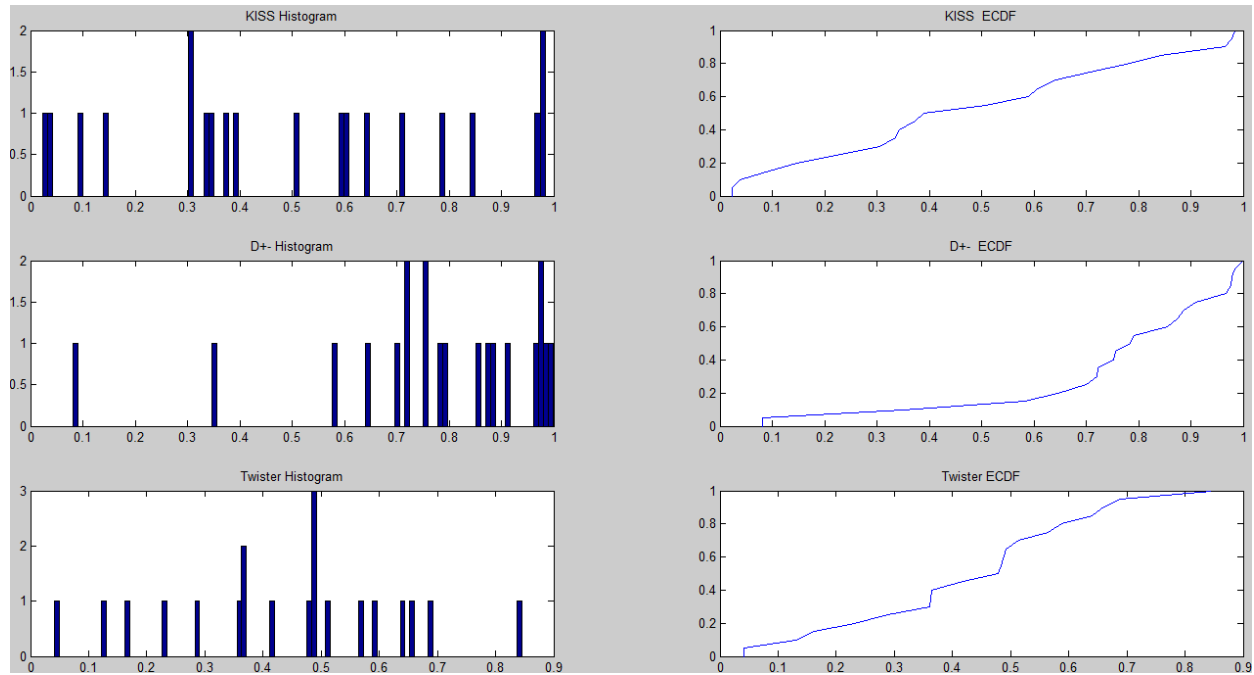

Figure 5

As you can see in Figure 5, the bitstream test does not appear to do very well as far as a uniform distribution for the D Sequence PRNG due to the values being concentrated toward the higher end of the the (0,1) range. Again, this alone does not constitute a failure, but when taken in conjunction with the prior results, one must begin to wonder how well this algorithm is doing as a whole.

*DNA, OPSO, and OQSO*

The results of the DNA, OPSO, and OQSO tests can be lumped together as they all came out essentially the same and they are extremely similar tests. The histogram and ECDF are both shown in Fig. 6:

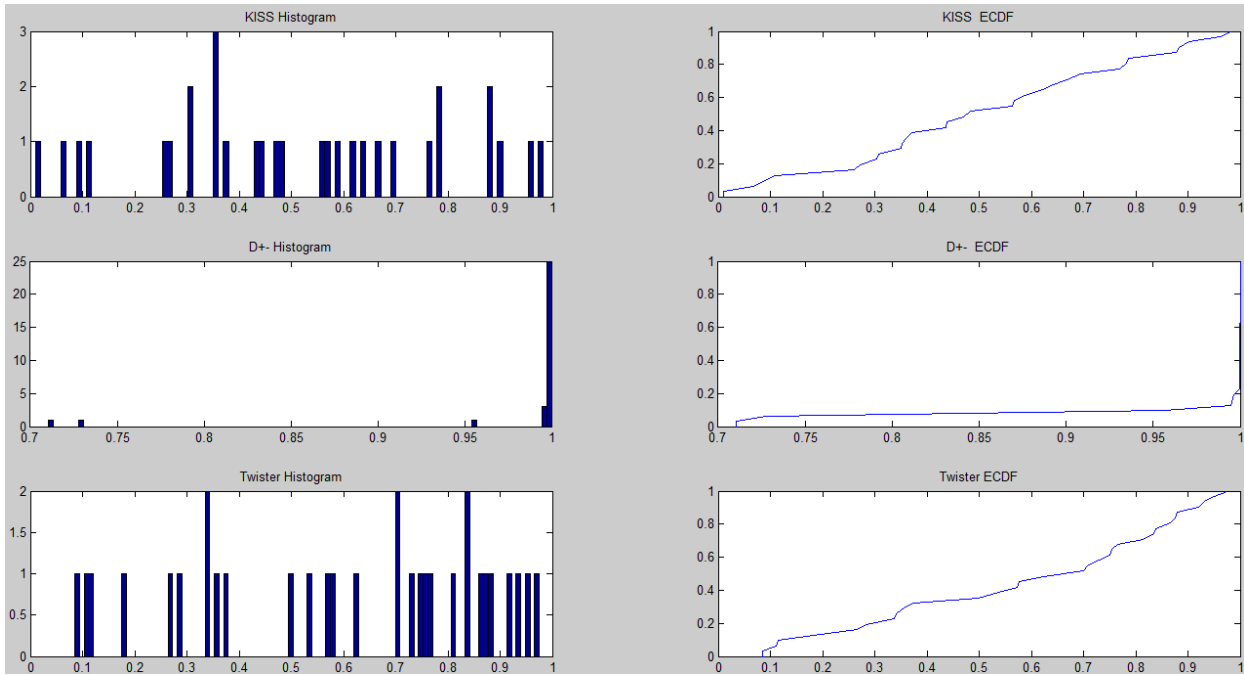

Figure 6

From this we can now see what a true failure looks like. While in the other tests there seemed to be some issues with values being close to the threshold of 'unacceptable' values, those were actually most likely just outliers that resulted from small sample size. However, in the case of these three tests I feel comfortable saying with 100% certainty that the D sequence PRNG failed.

*Squeeze Test*

| KISS | D Sequence | Twister |
| --- | --- | --- |
| .564635 | 1.000000 | .933079 |

Figure 7

Fig. 7 shows the results of the squeeze test, which would indicate a definite failure again coming from the D Sequence PRNG. This puts the total of failed tests at 2 for the D Sequence.

**Conclusions**

On the whole the D sequence algorithm did a fairly good job with these tests. Of the fifteen tests, I would only consider two failures. These two tests were the Squeeze Test and the DNA/OQSO/OPSO Test, as they were both fairly consistently 1.0000 values. Considering that the only other failures were single p-values that were not within six decimal places of 0 or 1 as mentioned by Marsaglia, I believe this is a reasonable assessment.